\documentclass[12pt,leqno]{article}

\usepackage{latexsym,amsmath,amssymb,amsthm,amsfonts}

\renewcommand{\epsilon}{\varepsilon}
\newtheorem{theorem}{Theorem}

\newtheorem{corollary}{Corollary}

\begin{document}

\author{Andriy Bondarenko\footnote{This work was carried out during the tenure of an ERCIM ``Alain Bensoussan'' Fellowship Programme. The research leading to these results has received funding from the European Union Seventh Framework Programme (FP7/2007-2013) under grant agreement $n\textsuperscript{o}$ 246016.}}
\title{On Borsuk's conjecture for two-distance sets}
\date{}
\maketitle
\begin{abstract}
In this paper we answer Larman's question on Borsuk's conjecture for two-distance sets.
We find a two-distance set consisting of $416$ points on the unit sphere $S^{64}\subset\mathbb{R}^{65}$ which cannot be partitioned into $83$ parts of smaller diameter. This also reduces the smallest dimension in which Borsuk's conjecture is known to be false. Other examples of two-distance sets
with large Borsuk's numbers will be given.
\end{abstract}
\section{Introduction}
\label{intro}
For each $n\in{\mathbb N}$ the Borsuk number $b(n)$ is the minimal number such that any bounded set in
$\mathbb{R}^{n}$ consisting of at least 2 points can be partitioned into $b(n)$ parts of smaller diameter.
In 1933 Karol Borsuk~\cite{Bor} conjectured that $b(n)=n+1$. The conjecture was disproved by Kahn and Kalai~\cite{KK}
who showed that in fact $b(n)>1.2^{\sqrt{n}}$ for large $n$. In particular, their construction implies that $b(n)>n+1$ for $n=1325$ and for all
$n>2014$. This result attracted a substantial amount of attention from many mathematicians; see for example~\cite{A},~\cite{BMP}, and~\cite{R}.
Improvements on the smallest dimension $n$ such that $b(n)>n+1$ were obtained by Nilli~\cite{N} ($n = 946$), Raigorodskii~\cite{R1} ($n = 561$), Wei\ss bach~\cite{W} ($n = 560$), Hinrichs~\cite{H} ($n = 323$), and Pikhurko~\cite{Pikh} ($n = 321$). Currently the best known result is that Borsuk's conjecture is false for $n\ge 298$; see~\cite{HR}. On the other hand, many related problems are still unsolved. Borsuk's conjecture can be wrong even in dimension $4$. Only the estimate $b(4)\le 9$ is known; see~\cite{La}.

In the 1970s Larman asked if the Borsuk's conjecture is true for two-distance sets; see also~\cite{K} and~\cite{R}.
Denote by $b_2(n)$ the Borsuk number for two-distance sets in dimension $n$, that is the minimal number such that any two-distance set in $\mathbb{R}^{n}$ can be partitioned into $b_2(n)$ parts of smaller diameter. The aim of this paper is to construct two-distance sets with large Borsuk's numbers. Two basic constructions follow from Euclidean representations of $G_2(4)$ and $Fi_{23}$ strongly regular graphs. First we prove
\begin{theorem}
There is a two-distance subset $\{x_1,\ldots,x_{416}\}$ of the unit sphere $S^{64}\subset\mathbb{R}^{65}$
such that $\langle x_i,x_j \rangle=1/5$ or $-1/15$ for $i\neq j$ which cannot be partitioned into $83$ parts of smaller diameter.
\end{theorem}
Hence $b(65)\ge b_2(65)\ge 84$. We also prove the following
\begin{theorem}
There is a two-distance subset $\{x_1,\ldots,x_{31671}\}$ of the unit sphere $S^{781}$
such that $\langle x_i,x_j\rangle =1/10$ or $-1/80$ for $i\neq j$ which cannot be partitioned into $1376$ parts of smaller diameter.
\end{theorem}
Then, using the configurations from Theorem 1 and Theorem 2 we prove
\begin{corollary}
For integers $n\ge 1$ and $k\ge 0$ we have
\begin{equation}
\label{car1}
b_2(66n+k)\ge 84n+k+1,
\end{equation}
and
\begin{equation}
\label{car2}
b_2(783n+k)\ge 1377n+k+1.
\end{equation}
\end{corollary}
Finally, using again the configuration from Theorem~2 we prove slightly better estimates for $b_2(781)$, $b_2(780)$, and $b_2(779)$ than what can be obtained by~\eqref{car1}.
\begin{corollary}
The following inequalities hold:
$$
b_2(781)\ge 1225,\quad b_2(780)\ge 1102,\quad\text{and}\quad b_2(779)\ge 1002.
$$
\end{corollary}
The paper is organized as follows. First, in Section~\ref{sec:1} we describe Euclidean representations of a strongly regular graph by
two-distance sets and then in Section~\ref{sec:2} we prove our main results.
\section{Euclidean representations of strongly regular graphs}
\label{sec:1}
A strongly regular graph $\Gamma$ with parameters
$(v,k,\lambda,\mu)$ is an undirected regular graph on $v$ vertices
of valency $k$ such that each pair of adjacent vertices
has $\lambda$ common neighbors, and each pair of nonadjacent
vertices has $\mu$ common neighbors. The adjacency matrix $A$ of
$\Gamma$ has the following properties:
$$
AJ = kJ
$$
and
$$
A^2 + (\mu - \lambda)A + (\mu - k)I = \mu J,
$$
where $I$ is the identity matrix and $J$ is the matrix with all
entries equal to~$1$ of appropriate sizes. These conditions imply
that
\begin{equation}
\label{par} (v - k - 1)\mu = k(k - \lambda - 1).
\end{equation}
Moreover, the matrix $A$ has only 3 eigenvalues: $k$ of multiplicity
$1$, one positive eigenvalue
$$
r=\frac 12\left(\lambda-\mu+\sqrt{(\lambda-\mu)^2+4(k-\mu)}\right)
$$
of multiplicity
\begin{equation}
\label{f}
f=\frac 12
\left(v-1-\frac{2k+(v-1)(\lambda-\mu)}{\sqrt{(\lambda-\mu)^2+4(k-\mu)}}\right),
\end{equation}
and one negative eigenvalue
$$
s=\frac 12\left(\lambda-\mu-\sqrt{(\lambda-\mu)^2+4(k-\mu)}\right)
$$
of multiplicity
\begin{equation*}
\label{g} g=\frac 12
\left(v-1+\frac{2k+(v-1)(\lambda-\mu)}{\sqrt{(\lambda-\mu)^2+4(k-\mu)}}\right).
\end{equation*}
Clearly, both $f$ and $g$ must be integers. This together
with~\eqref{par} gives a collection of feasible parameters
$(v,k,\lambda,\mu)$ for strongly regular graphs.

Let $V$ be the set of vertices $\Gamma$. Consider the columns $\{y_i:
i\in V\}$ of the matrix $A-sI$ and put $x_i:=z_i/\|z_i\|$, where
$$
z_i=y_i-\frac 1{v}\sum_{j\in V}y_j, \quad i\in V.
$$
Note that while the vectors $x_i$ lie in $\mathbb{R}^v$, they span
at most an $f$-dimensional vector space. Thus for convenience we
consider them to lie in $\mathbb{R}^f$. By easy calculations
$$
\langle x_i,x_j\rangle=
\begin{cases}1, & \mbox{if }i=j, \\
p, & \mbox{if }i\mbox{ and }j\mbox{ are adjacent},\\
q, & \mbox{otherwise},
\end{cases}
$$
where
\begin{equation}
\label{pq}
p=\frac{\lambda-2s-\beta}{s^2+k-\beta},\quad q=\frac{\mu-\beta}{s^2+k-\beta},\quad\beta=\frac 1v(s^2+k+k(\lambda-2s)+(v-k-1)\mu).
\end{equation}
Denote by $\Gamma_f$ the configuration $x_i$, $i\in V$. Similarly, we can define the configuration $\Gamma_g$ in $\mathbb{R}^g$. The configurations $\Gamma_f$ and $\Gamma_g$ were also considered in~\cite{Cam} and have many other fascinating properties. For example, they are spherical $2$-designs.
\section{Proof of main results}
\label{sec:2}
For any vertex $v\in V$ of a strongly regular graph $\Gamma$, let $N(v)$ be the set of all neighbors of $v$
and let $N'(v)$ be the set of non-neighbors of $v$, i.e.
$N'(v)=V\setminus(\{v\}\cup N(v))$.\\

\noindent
{\em Proof of Theorem 1}
We consider the configuration $\Gamma_f$ of the well-known strongly regular graph $\Gamma=G_2(4)$ with parameters $(416,100,36,20)$. By~\eqref{f} we have that $f=65$. Moreover by~\eqref{pq}, $p=1/5$ and $q=-1/15$. Therefore the diameter of $\Gamma_f$ is the distance between $x_i$ and $x_j$ where $i$ and $j$ are nonadjacent. Hence, the configuration cannot be partitioned into less than $v/m$ parts, where $m$ is the size of the largest clique in $\Gamma$.
To prove Theorem 1 it is enough to show that $\Gamma$ has no $6$-clique.
We will use the following theorem consisting of four independent results that can be found in~\cite{B}.\\

\noindent
{\bf Theorem A}
{\it \begin{itemize}
\item[(i)] For each $u\in V$ the subgraph of $\Gamma$
induced on $N(u)$ is a strongly regular graph with parameters $(100,36,14,12)$ (the Hall-Janko graph).
In other words the Hall-Janko graph is the first subconstituent of $\Gamma$.
 \item[(ii)]The first subconstituent of the Hall-Janko graph is the $U_3(3)$ strongly regular graph with parameters $(36,14,4,6)$.
\item[(iii)] The first subconstituent of $U_3(3)$ is a graph on $14$ vertices of regularity $4$ (the co-Heawood graph).
\item[(iv)] The co-Heawood graph has no triangles.
\end{itemize}
}
Parts (i)-(iii) are folklore. They follow from D.G. Higman's theory of rank 3 permutation groups (see also~\cite{G} and~\cite{L}).
Part (iv) follows from the fact that the co-Heawood graph is a subgraph of the Gewirtz graph with parameters (56,10,0,2); see also~\cite{BCD}.

Now, for vertices $u,v,w\in V$ forming a triangle, (i)-(iii) imply that $$|N(u)\cap N(v)\cap N(w)|=14.$$ Moreover, the subgraph induced on $N(u)\cap N(v)\cap N(w)$ is the co-Heawood graph. Therefore by (iv) the maximal cliques in $\Gamma$ are of size 5.
\qed

\noindent
{\em Proof of Theorem 2}
Consider the configuration $\Gamma_f$ of the $Fi_{23}$ graph with parameters $(31671,3510,693,351)$.
We have $f=782$, $p=1/10$, and $q=-1/80$. Hence, the diameter of $\Gamma_f$ is the distance between nonadjacent vertices.
Therefore $\Gamma_f$ cannot be partitioned into less than $v/m$ parts, where $m$ is the size of the largest clique in $\Gamma$. We will use the well-known fact (see~\cite{P}) that the first subconstituent of $\Gamma$ is the strongly regular graph with parameters $(3510,693,180,126)$ and the second subconstituent of $\Gamma$ is the strongly regular graph $G$ with parameters $(693,180,51,45)$. Now we will estimate from above the size of a clique in $G$. To this end consider the complement graph $\bar G$ having parameters $(693,512,376,384)$. For the configuration $\bar G_f$, we have that $f=440$, $p=1/64$, and $q=-1/20$. Therefore, the size of a clique $K$ in $G$ cannot be larger than $21$. Otherwise the vector
$$\sum_{i\in K}x_i,\qquad x_i\in\bar G_f,$$
is of negative norm. Thus, the size of a clique in $\Gamma$ is not larger than $23$ and hence $\Gamma_f$ cannot be partitioned into less than $31671/23=1377$ parts of smaller diameter.
\qed

\noindent
{\em Proof of Corollary 1}
Let us first prove~\eqref{car1} for $k=0$. Fix $n\in\mathbb{N}$ and put $m=66n$.
Consider the following coordinate representation of a vector $y\in\mathbb{R}^{m}$:
$$
y=(y_1,\ldots,y_n|a_1,\ldots,a_n),
$$
where $y_k\in\mathbb{R}^{65}$ and $a_k\in\mathbb{R}$, $k=1,\ldots, n$.
Now we take the following set of unit vectors in $\mathbb{R}^{m}$:
$Y=\{v_{ik},\,i=1,\ldots, 416,\,k=1,\ldots, n\}$, where
$$
v_{ik}=(0,\ldots,0,\frac{\sqrt{15}}4x_i,0,\ldots,0\,|\,0,\ldots,0,\frac 14,0,\ldots,0),\, i=1,\ldots,416,\,k=1,\ldots, n,
$$
Here each $v_{ik}$ has only two nonzero coordinates $y_k$ and $a_k$, and vectors $x_i$ are such as in Theorem 1.
Clearly, $\langle v_{ik},v_{jl}\rangle=0$ if $k\neq l$. Moreover,
$$
\langle v_{ik},v_{jk}\rangle=
\begin{cases}1, & \mbox{if }i=j, \\
1/4, & \mbox{if }i\mbox{ and }j\mbox{ are adjacent},\\
0, & \mbox{otherwise}.
\end{cases}
$$
Therefore, $Y$ is a two-distance set consisting of $416n$ vectors.
Now, by Theorem~1, this set cannot be partitioned into less than $84n$ parts of smaller diameter. Adding the vector
$v$ which is at distance $\sqrt{2}$ to each vector of~$Y$
$$
v=(0,\ldots,0\,|\,\alpha,\ldots,\alpha),\quad \alpha=\frac{1+\sqrt{1+16n}}{4n}
$$
($\alpha$ is a solution of the equation $(\alpha-1/4)^2+(n-1)\alpha^2=17/16$) we obtain that $b_2(m)\ge 84n+1$.
Finally we note that all these $416n+1$ vectors are at the same distance $R$ to the vector $(0,\ldots,0\,|\,\gamma,\ldots,\gamma)$,
where
$$
\gamma=\frac{\alpha}{4n\alpha-1}\,\text{ and }\, R=\frac{4\sqrt{n}}{\sqrt{16n+1}}<1
$$
($\gamma$ is a solution of the equation $(\gamma-1/4)^2+(n-1)\gamma^2+15/16=n(\alpha-\gamma)^2$).
Hence we can add a new vector at the diameter distance $\sqrt{2}$ to each of these $416n+1$ vectors
to get a new set of $416n+2$ vectors in $\mathbb{R}^{m+1}$ provided that $b_2(m+1)\ge 84n+2$.
We can also rescale this new set to be on the sphere $S^{m}$. Now inductive application of this procedure immediately gives us~\eqref{car1}. This procedure was also described in~\cite[Lemma 9]{HR}.
Similarly, Theorem 2 implies~\eqref{car2}.
\qed

\noindent
{\em Proof of Corollary 2}
Let $\Gamma$ be the $Fi_{23}$ graph. For a vertex $u\in V$, consider the subset $\left
\{x_i:\, i\in N'(u)\right\}$ of the configuration $\Gamma_f$.
This subset lies in the hyperplane $\langle x_u,x\rangle=-1/80$ and consists of $31671-3510-1=28160$ vectors. Hence,
$b_2(781)>[28160/23]=1224$.

Similarly, for adjacent vertices $u$ and $v$, the subset $\left
\{x_i:\, i\in N'(u)\cap N'(v)\right\}$ consists of $31671-2\times3510+693=25344$ vectors. This subset lies
in the hyperplane $\{x\in\mathbb{R}^{782}:\,\langle x_u,x\rangle=-1/80\text{ and }\langle x_v,x\rangle=-1/80\}$, and hence $b_2(780)>[25344/23]=1101$.

Finally, consider a subset $\left
\{x_i:\,i\in N'(u)\cap N'(v)\cap N'(w)\right\}$ such that the vertices $u$, $v$, $w$ form a triangle. This subset consists of $31671-3\times3510+3\times693-180=23040$ vectors, and hence $b_2(779)>[23040/23]=1001$.
\qed\\
{\bf Acknowledgements}
The author thanks Danylo Radchenko and Kristian Seip for several fruitful discussions.
The author also wish to thank the Centre for Advanced Study at the Norwegian Academy of Science and Letters in Oslo, and Mathematisches Forschungsinstitut Oberwolfach for their hospitality during the preparation of this manuscript and for providing a stimulating
atmosphere for research.

\vspace{0.2cm}
\noindent Andriy V. Bondarenko\\
Department of Mathematical Analysis,\\ Taras Shevchenko National
University of Kyiv, str.\ Volodymyrska, 64,\\ Kyiv, 01033, Ukraine\\
and\\
Department of Mathematical Sciences, Norwegian University of Science and Technology, NO-7491 Trondheim, Norway\\
{\it Email address: andriybond@gmail.com}

\end{document}